\documentclass[review,onefignum,onetabnum]{siamart171218}

\usepackage{lipsum}
\usepackage{amsfonts}
\usepackage{graphicx}
\usepackage{epstopdf}
\usepackage{algorithmic}
\ifpdf
  \DeclareGraphicsExtensions{.eps,.pdf,.png,.jpg}
\else
  \DeclareGraphicsExtensions{.eps}
\fi
\usepackage{subfig}

\newsiamremark{remark}{Remark}
\newsiamremark{hypothesis}{Hypothesis}
\crefname{hypothesis}{Hypothesis}{Hypotheses}
\newsiamthm{claim}{Claim}

\headers{Taylor Polynomials in High Precision as Universal Approximators}{Nikolaos P. Bakas}

\title{Taylor Polynomials in High Arithmetic Precision as Universal Approximators\thanks{Submitted to the editors DATE.
}}

\author{Nikolaos P. Bakas\thanks{Intelligent Systems Lab \& Civil Engineering Department, Neapolis University Pafos, Danais 2, 8042, Pafos, Cyprus 
  (\email{n.bakas@nup.ac.cy}, \url{http://www.nup.ac.cy/}).}}

\usepackage{amsopn}

\makeatletter
\newcommand*{\addFileDependency}[1]{
  \typeout{(#1)}
  \@addtofilelist{#1}
  \IfFileExists{#1}{}{\typeout{No file #1.}}
}
\makeatother

\newcommand*{\myexternaldocument}[1]{%
    \externaldocument{#1}%
    \addFileDependency{#1.tex}%
    \addFileDependency{#1.aux}%
}

\nolinenumbers

\ifpdf
\hypersetup{
  pdftitle={Taylor Polynomials in High Arithmetic Precision as Universal Approximators},
  pdfauthor={Nikolaos P. Bakas}
}
\fi

\myexternaldocument{ex_supplement}

\usepackage{dirtytalk}

\begin{document}

\maketitle

\begin{abstract}
Function approximation is a generic process in a variety of computational problems, from data interpolation to the solution of differential equations and inverse problems. In this work, a unified approach for such techniques is demonstrated, by utilizing partial sums of Taylor series in high arithmetic precision. In particular, the proposed method is capable of interpolation, extrapolation, numerical differentiation, numerical integration, solution of ordinary and partial differential equations, and system identification. The method is based on the utilization of Taylor polynomials, by exploiting some hundreds of computer digits, resulting in highly accurate calculations. Interestingly, some well-known problems were found to reason by calculations accuracy, and not methodological inefficiencies, as supposed. In particular, the approximation errors are precisely predictable, the Runge phenomenon is eliminated and the extrapolation extent may a-priory be anticipated. The attained polynomials offer a precise representation of the unknown system as well as its radius of convergence, which provide a rigor estimation of the prediction ability. The approximation errors have comprehensively been analyzed, for a variety of calculation digits and test problems. 
\end{abstract}



\section{Introduction}
\label{Introduction}
The utilization of High Arithmetic Precision (HAP) for the modeling of an unknown function exhibited a remarkable extrapolation ability in \cite{Bakas2019}, with extrapolation spans of 1000\% higher than the existing methods in the literature. The origin of this method was the modeling of an unknown analytic function, which is an essential issue in a variety of numerical methods, with high arithmetic precision. Standard programming languages are limited to 16 to 64 floating point digits, and researchers have been taking into account high arithmetic precision for the various computations regarding numerical integration \cite{Bailey2005}, interpolation \cite{Cheng2012} and solution of Partial Differential Equations (PDEs) \cite{Huang2007}, however high arithmetic precision has not been studied extensively yet. To the contrary, standard techniques exist for interpolation with Taylor polynomials \cite{Guessab2006162,Kalantari2000287}, as well as the solution of differential equation \cite{Berz1998361,Yalçinbaş2000291}, however, certain problems occur, as the well-known Runge phenomenon \cite{Platte2011,Boyd199257}, which remains a major complication \cite{Zhang20161948,Boyd2009158,Boyd2009484}. 

 Taylor series arise in the foundation of Differential Calculus \cite{Taylor1715}, by associating the behavior a function around a point $x_0$, with its derivatives on that particular point. Despite the vast literature on the function approximation with Taylor series as well as their instabilities, no analysis and discussion exist on their theoretical explanation. Accordingly, although Taylor series are capable of approximating any analytic function, because in practice they often fail, and researchers use other approximators than Taylor polynomials, such as Radial Basis Functions, Lagrange Polynomials, Chebyshev Polynomials, Artificial Neural Networks, etc, to avoid numerical instabilities. A variety of Numerical Methods have been developed for such operations, as researchers have been observing that Taylor polynomials do not offer stable calculations. Utilizing high-arithmetic precision, we demonstrate that such need, which arose to cover the computational inaccuracies, does not exist. Taking into account the high extrapolation spans attained in \cite{Bakas2019}, obtained with integrated radial basis functions \cite{Babouskos2015,Yiotis2015} and some hundreds or even thousands of digits for the calculations, we applied high arithmetic precision, utilizing the BigFloat structure of Julia Language \cite{bezanson2017julia}, to truncated Taylor series, known as Taylor Polynomials or Partial Sums.

The purpose of this work was to present a unified approach for the interpolation, extrapolation, numerical differentiation, solution of partial differential equations, system identification and numerical integration for problems which supply only some given data of the unknown analytic function or the source for PDEs. The paper is organized as follows. The formulation of our approach is presented in \cref{Description-method}, some basic operations and results for !-Dimensional Interpolation, Extrapolation, Numerical Differentiation, Numerical Integration, solution of Ordinary Differential Equations are in \cref{Function-approximation}, results for multidimensional Function Approximation, solution of Partial Differential Equations and System Identification are in \cref{N-D}, and the conclusions follow in
\cref{sec:conclusions}.

\section{Description of the method}
\label{Description-method}

Let $f(x)$ be an analytic function, which is unknown. It is given that the function takes values $\mathbf{f}=\left\{ {{f}_{1}},{{f}_{2}},...,{{f}_{N}} \right\}$ at specified points $\mathbf{x}=\left\{ {{x}_{1}},{{x}_{2}},...,{{x}_{N}} \right\}$ as in Figure 1, for a generic analytic function. By applying the Taylor series \cite{Taylor1715,wikiTaylor} of the function at some point ${{x}_{0}}$, we may write $f(x\pm {{x}_{0}})=f({{x}_{0}})\pm \frac{{f}'({{x}_{0}})}{1!}(x-{{x}_{0}})+\frac{{f}''({{x}_{0}})}{2!}{{(x-{{x}_{0}})}^{2}}\pm \cdots \pm \frac{{{f}^{(n)}}({{x}_{0}})}{n!}{{(x-{{x}_{0}})}^{n}}+\cdots$. The derivatives of the function, $\mathbf{df}=\left\{ {{f}^{0}},{f}',{f}'',...,{{f}^{(n)}} \right\}$ at ${{x}_{0}}$, divided by $n!$, are constant quantities, hence by truncating the series at the ${{n}^{th}}$ power, we derive that
 \[f(x\pm {{x}_{0}})\cong f({{x}_{0}})\pm \frac{{f}'({{x}_{0}})}{1!}(x-{{x}_{0}})+\frac{{f}''({{x}_{0}})}{2!}{{(x-{{x}_{0}})}^{2}}\pm \cdots +\frac{{{f}^{(n)}}({{x}_{0}})}{n!}{{(x-{{x}_{0}})}^{n}}+{{R}_{n}}(x).\]
 The remainder of the approximation is bounded \cite{apostol1967calculus,wikiRemaind} by \[|{{R}_{n}}(x)|\le \frac{{{f}^{n+1}}(x)}{(n+1)!}|x-{{x}_{0}}{{|}^{n+1}},\forall x:|x-{{x}_{0}}|\le r.\]\par 
 
 For a series$f(x)=\sum_{n=0}^{\infty }{{{a}_{n}}}{{(x-{{x}_{0}})}^{n}}$, we have that the radius of convergence \cite{apostol1967calculus,wikiRadius} $r$, is a non-negative real number or $\infty $ such that the series converges if $|x-{{x}_{0}}|<r$ and diverges if $|x-{{x}_{0}}|\ge r$, that is to say, the series converges in the interval $({{x}_{0}}-r,{{x}_{0}}+r)$. We may compute $r$ by the ratio test $\lim \sup \left| {{a}_{n+1}}/{{a}_{n}} \right|$ or by the root test, with 
 $r=1/\limsup_{n \to \infty} \sqrt[n]{|{{a}_{n}}|}.$
 We select the root test because the coefficients $a_i$ many times contain zero elements and the division is not computationally stable. Furthermore, because $\lim \inf ({{a}_{n+1}}/{{a}_{n}})\le \lim \inf ({{({{a}_{n}})}^{(1/n)}})\le \lim \sup ({{({{a}_{n}})}^{(1/n)}})\le \lim \sup ({{a}_{n+1}}/{{a}_{n}})$ \cite{browder2012mathematical}, the computes $r$ are from the root test higher than the ratio test. High arithmetic precision, found capable for the accurate computation of $r$, for known series, while floating-point fails. This is a significant part of the proposed numerical schemes as, the identification of $r$, offers information on the larger disk where the series converges. Accordingly, we obtain knowledge of the interpolation accuracy or even the extrapolation span of the approximated function beyond the given domain. \par 
 
In particular, at ${{x}_{0}}=0$, we may write that
\begin{equation}f(x)\cong {{a}_{0}}\pm {{a}_{1}}x+{{a}_{2}}{{x}^{2}}\pm \cdots +{{a}_{n}}{{x}^{n}}\end{equation} where $\mathbf{a}=\left\{ 1,{f}'/1,{f}''/2!,...,{{f}^{(n)}}/n! \right\}=\mathbf{df}./\{1,\ldots ,n!\}$. This is the truncated Taylor polynomial, which may converge to $f$ \cite{Katsoprinakis2011,Nestoridis2011}. By applying the Taylor formula for all the $n$ given points ${{x}_{i}}$, with $i=1\ldots n$, we obtain
$\mathbf{f}=\mathbf{Va}$
where $\mathbf{V}$is the Vandermonde matrix, with elements $v_{i,j}={x_i}^{j-1}$, where $j=1\ldots n$ \cite{press2007vwt,horn1991topics,wikiVand}.

\begin{figure}[ht]
        \begin{tabular}{p{5cm}c}
            {$\mathbf V = 
                {\renewcommand{\arraystretch}{1.2}
                \begin{pmatrix}
           1 & {{x}_{1}} & x_{1}^{2} & \ldots  & x_{1}^{n-1}  \\
           1 & {{x}_{2}} & x_{2}^{2} & \ldots  & x_{2}^{n-1}  \\
           1 & {{x}_{3}} & x_{3}^{2} & \ldots  & x_{3}^{n-1}  \\
           \vdots  & \vdots  & \vdots  & \ddots  & \vdots   \\
           1 & {{x}_{n}} & x_{n}^{2} & \ldots  & x_{n}^{n-1}  \\
        \end{pmatrix}} 
            $}
            &
            $\vcenter{\hbox{\includegraphics[width=50mm,scale=0.5]{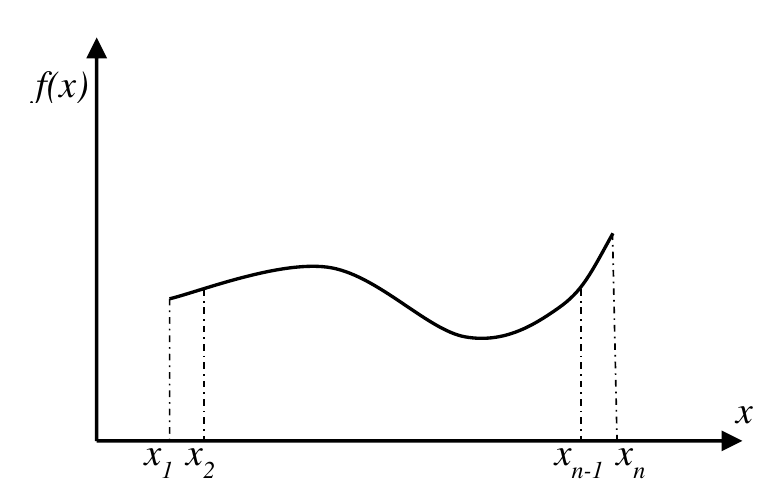}}}$
        \end{tabular}
        \caption{Given values of $f(x)$ at points $x_i$ for the approximation of $f$ by inverting the corresponding Vandermonde matrix $\mathbf V$.}
        \label{Figure 1}
\end{figure}

The square Vandermonde matrix for distinct $x_i$ is invertible, with $\det (\mathbf{V})=\prod_{1\le i<j\le n}{({{x}_{j}}-{{x}_{i}})}$ \cite{ycart2012case}, and inverse matrix  ${{\mathbf{V}}^{\mathbf{-1}}}\mathbf{=}{{\mathbf{U}}^{\mathbf{-1}}}{{\mathbf{L}}^{\mathbf{-1}}}$, where the elements ${{l}_{ij}}$ of the ${{\mathbf{L}}^{\mathbf{-1}}}$, and ${{u}_{ij}}$ of ${{\mathbf{U}}^{\mathbf{-1}}}$, are given by
${{l}_{ij}}=\left\{ \prod_{k=1(k\ne j)}^{i}{\frac{1}{{{x}_{j}}-{{x}_{k}}}};0\forall i<j;{{l}_{11}}=1 \right\}$,${{u}_{ij}}=\left\{ {{u}_{i-1,j-1}}-{{u}_{i,j-1}}{{x}_{j-1}};{{u}_{i1}}=0;{{u}_{ii}}=1,{{u}_{oj}}=0 \right\}$ \cite{turner1966inverse}.
Hence we have closed-form formulas for the matrix ${{\mathbf{V}}^{\mathbf{-1}}}$, and for $\det (\mathbf{V})$, which will be later used for the comparison among the various digits utilized in the calculations. Accordingly, we can compute the polynomial factors $\mathbf{a}=\left\{ {{a}_{1}},{{a}_{2}},...,{{a}_{n}} \right\}$, by \[\mathbf{a}={{\mathbf{V}}^{-1}}\mathbf{f}.\] The computation of $\mathbf{a}$ with floating-point arithmetic exhibits significant errors in the inversion as well as the determinant calculation, with respect to their theoretical values by the closed-form formulas and numerical values computed by the computer.

\section{Function approximation in HAP}
\label{Function-approximation}
We will demonstrate the proposed numerical scheme, in a variety of numerical methods, analytic functions, and calculation digits. We begin with some basic operations.

\subsection{Basic Operations}
For the simple function$f(x)=\sin (x)$, the theoretical Taylor series exhibits alternating sign with intermediate zero coefficients \[\sin{x}=\sum\limits_{n=0}^{\infty }{\frac{{{(-1)}^{n}}}{(2n+1)!}}{{x}^{2n+1}}=0+x+0-\frac{{{x}^{3}}}{3!}+0+\frac{{{x}^{5}}}{5!}-\ldots ,\] hence according to the presented method the factors $\mathbf{a}=\left\{ {{a}_{1}},{{a}_{2}},...,{{a}_{n}} \right\}$, should be equal to $\left\{0,1,0,-\frac{1}{3!},0,\frac{1}{5!},-\ldots, \frac{1}{n!} \right\}$, for a truncated series with $n$ terms. However, the computation of ${{\mathbf{V}}^{\mathbf{-1}}}$, as well as the $\det (\mathbf{V})$, exhibits great variation with the calculation precision in bits $p$ (approximately equivalent to $p/3$ digits), when computed numerically or analytically by formulas. Table 1 presents such variation for $f(x)=\sin (x)$, with $L=1,n=201,dx=2L/(n-1)={{10}^{-2}}$, and $x\in [-L,L]$. The subscript \say{an} denotes the analytical value and \say{nu} the numerical one, as computed in variable precision $p=50 \ \text{to} \  2000$ bits.

\overfullrule=0pt
\begin{table}[htbp]
\caption{\label{tab:table-1}Variation of ${{\mathbf{V}}^{\mathbf{-1}}}$, $\det (\mathbf{V})$, and $\mathbf a$,with the calculation precision in bits $p$, for the same example.}
\begin{tabular}{ |p{2cm}||p{1.55cm}|p{1.5cm}|p{1.5cm}|p{1.5cm}|p{1.5cm}|  }
 \hline
 & $p=50$ & $p=100$ & $p=500$ & $p=1000$ & $p=2000$ \\
 \hline
$\det{{\mathbf{V}}_{an}} - \det{{\mathbf{V}}_{nu}}$ & 3.866e-2341	& 4.300e-4106 & -2.735e-6810 & -3.741e-6960 & -1.853e-7261\\
\hline
$\max | {{\mathbf{V}}^{\mathbf{-1}}}_{an}-{{\mathbf{V}}^{\mathbf{-1}}}_{nu} |$ & 9.739e+100 & 4.911e+94 & 1.242e+38 & 1.124e-111 & 5.504e-413\\
\hline
$\max | {{\mathbf{a}}_{an}}-{{\mathbf{a}}_{nu}} |$ & 4.029e+01 & 1.813e+00 & 9.252e-18 & 9.252e-18 & 9.252e-18\\
\hline
\end{tabular}
\end{table}

In Table 1, a high variation of the differences among ${{\mathbf{V}}^{\mathbf{-1}}}_{an}$and ${{\mathbf{V}}^{\mathbf{-1}}}_{nu}$ is revealed, from  9.739e+100 for $p=50$ bits, which is approximately equal to Floating-Point Precision, to 5.504e-413 for $p=2000$ bits. Accordingly, the maximum differences between $\mathbf a^{-1}_{an}$and $\mathbf a^{-1}_{nu}$ are 4.029e+01 for $p=50$ bits, and 9.252e-18 for $p>=500$ bits. It is important to underline that all the calculation regard the same example and same approximation scheme. Apparently, the errors of $O(\>10^{-16})$ cannot be considered as negligible. The significance of the precise computation is further demonstrated for the corresponding differences in the calculation of the determinant, with an analytical value constant at 1.647e-6754 and the corresponding differences from the computed, varying from 3.866e-2341 to -1.853e-7261, with alternating signs, again for the same example. In Table \ref{tab:table-1}, we also present that as the determinants' difference shortens, the same stands for the inversion errors.

\begin{figure}[ht]   
\centering
\subfloat[$p=2000$ bits]{\includegraphics[width=0.5\textwidth, keepaspectratio]{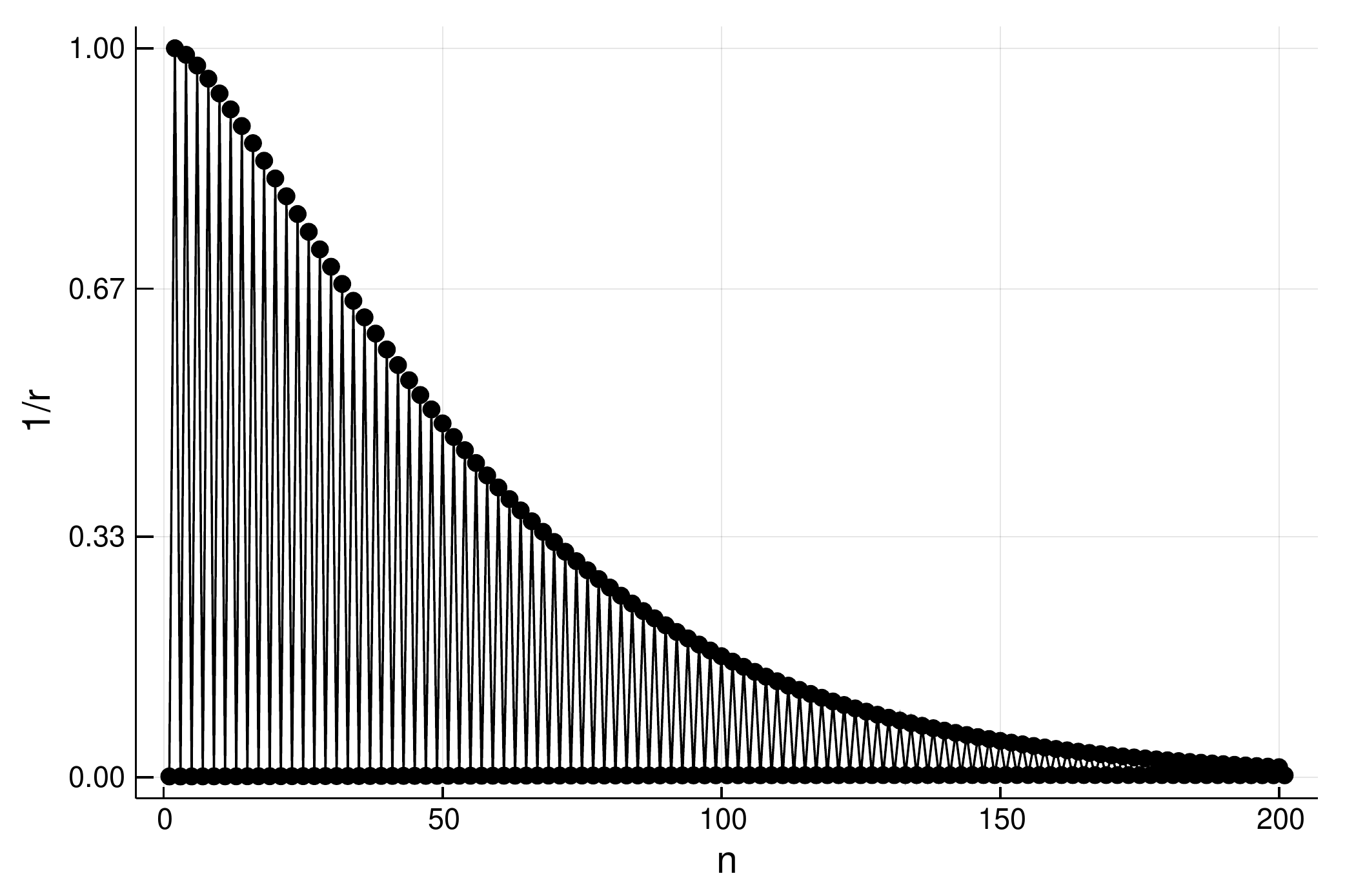}\label{fig:radius-high}}
\subfloat[$p=50$ bits]{\includegraphics[width=0.5\textwidth, keepaspectratio]{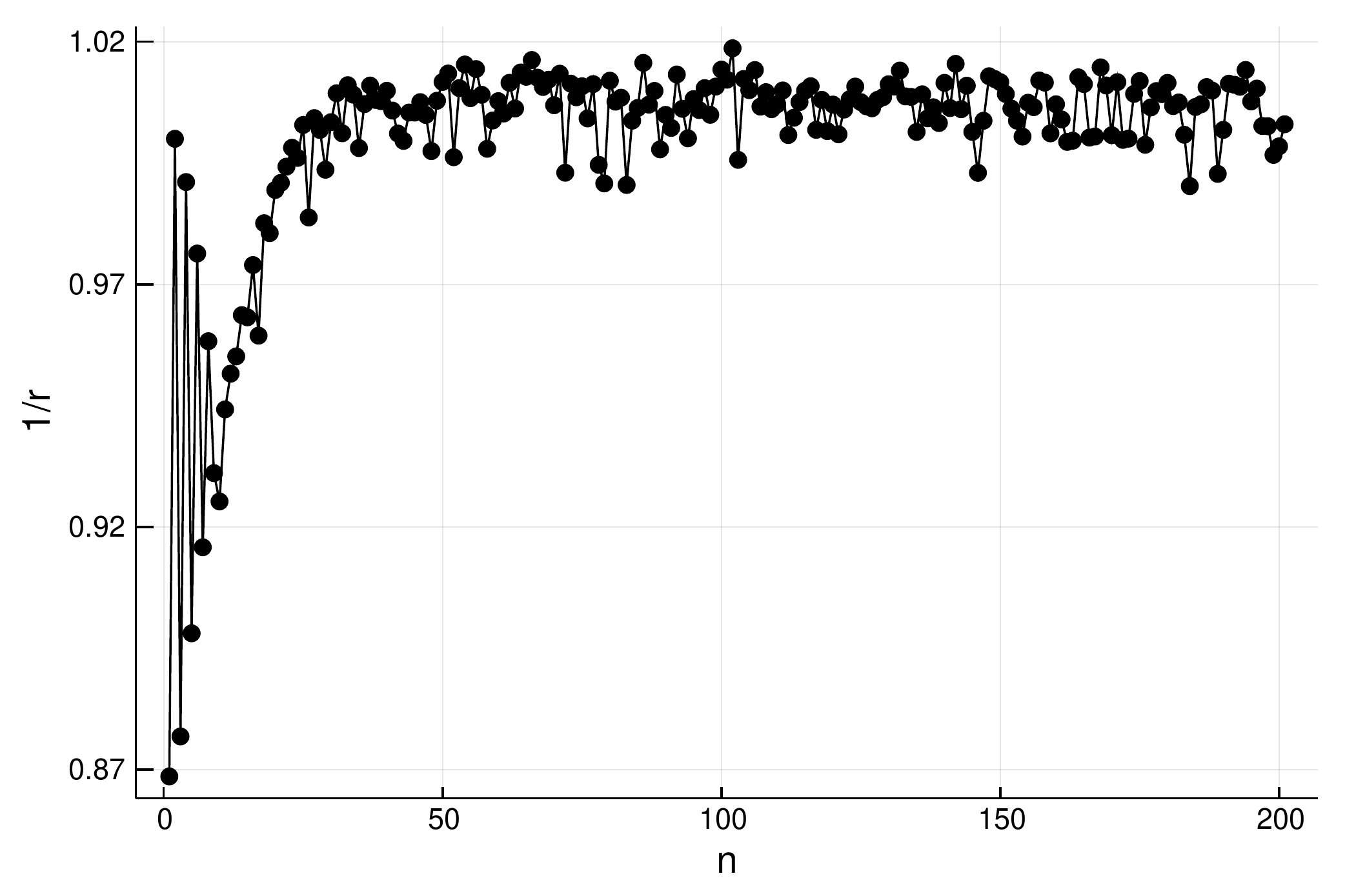}\label{fig:radius-low}}
\caption[Optional caption for list of figures 5-8]{Radius of convergence for the computed Taylor expansion of $f(x)$. }
\label{fig:radius}
\end{figure}

Digits accuracy exhibits great variation among the computed $1/r$ also. Precise calculation of $\mathbf{V}$and ${{\mathbf{V}}^{\mathbf{-1}}}$ makes convergent the computation of $1/r$, as the calculated ${\mathop{\lim \sup_{n\to \infty } }}\,\sqrt[n]{|{{a}_{n}}|}\simeq {\mathop{\lim \inf_{n\to \infty } }}\,\sqrt[n]{|{{a}_{n}}|}$ (Figure \ref{fig:radius-high}). Similarly, for the vector $\mathbf{a}$, the maximum absolute differences among analytical and numerical vary between 4.029e+01 and 9.252e-18. 

\subsection{Function Approximation}
As $f(x)=\sin (x)$, we have that $\left| {{f}^{n+1}}(x) \right|\le 1$], hence the theoretical remainder of the approximation, when using $n$terms of the Taylor series, is bounded as $|{{R}_{n}}(x)|\le \frac{1}{(n+1)!}|1-0{{|}^{n+1}}=6.308\text{e-378}$. In Table 2, the differences among computed and analytical values of $f$ at $x$ and ${{x}_{i}}=x+dx/2$ are presented.

\begin{table}[htbp]
\caption{\label{tab:table-2}Variation of approximation errors with the calculation precision in bits $p$.}
\begin{tabular}{ |p{2cm}||p{1.5cm}|p{1.5cm}|p{1.5cm}|p{1.5cm}|p{1.5cm}|  }
 \hline
 & $p=50$ & $p=100$ & $p=500$ & $p=1000$ & $p=2000$ \\
 \hline
$\max | f(x)_{an}-f(x)_{nu} |$ & 1.708e-12 & 3.045e-28 & 1.231e-148 & 3.770e-299 & 3.475e-600\\
\hline
$\max | f(x_i)_{an}-f(x_i)_{nu} |$ & 5.932e-08 & 2.045e-15 & 3.673e-96 & 2.373e-246 & 9.909e-407\\
\hline
\end{tabular}
\end{table}

Interestingly, although for $p=50$, the approximation error for $f(x)$ on the given points $x$, is 1.708e-12, the corresponding interpolation error on ${{x}_{i}}$, is 5.932e-08 (Table 2). The Runge phenomenon, which is severe at the boundaries, is eliminated, for $p>500$.

\begin{figure}[ht]   
\centering
\includegraphics[width=0.7\textwidth, keepaspectratio]{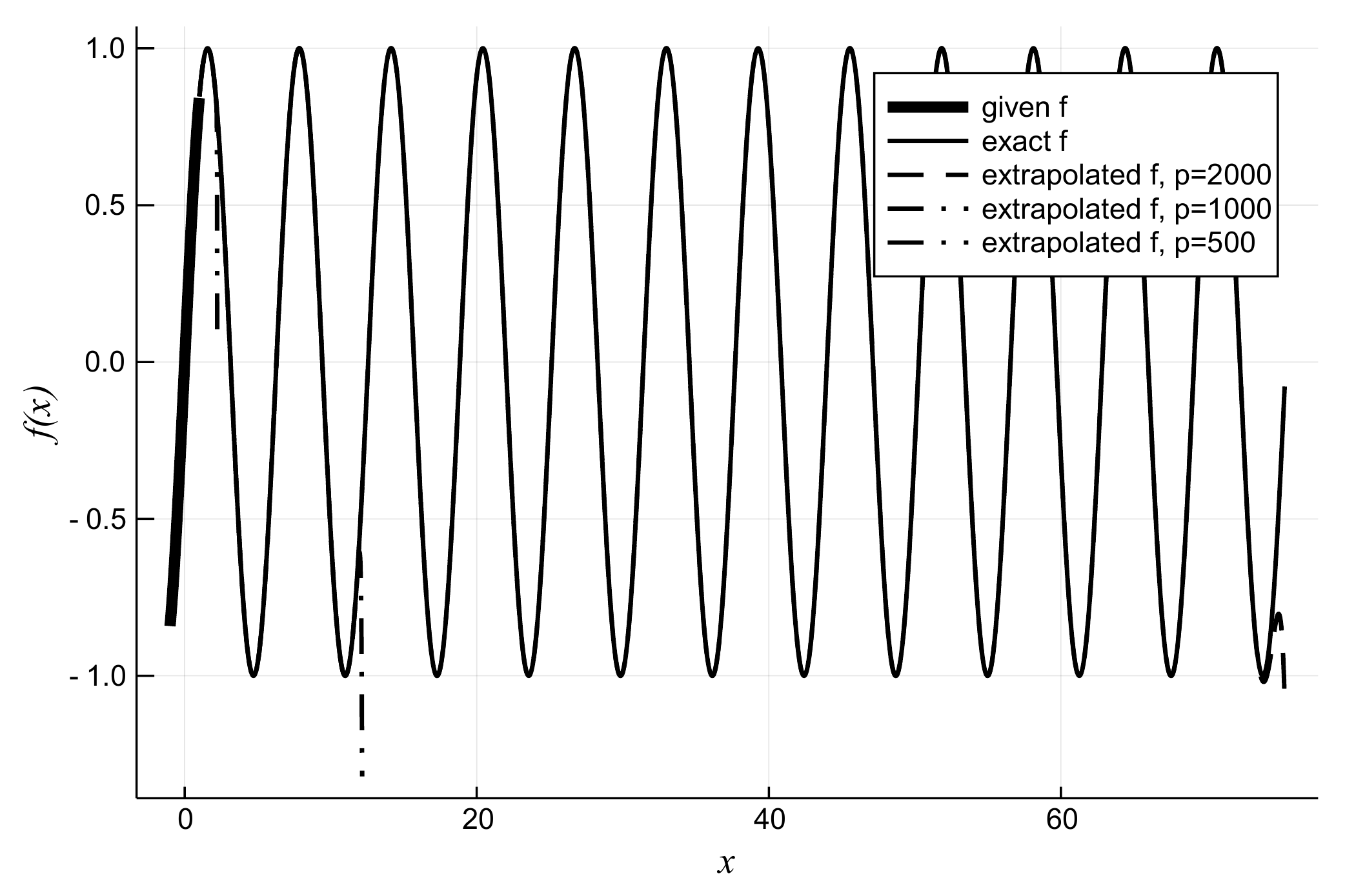}
\caption[]{Extrapolation of $f$, for varying arithmetic precision $p$. }
\label{fig:extrap}
\end{figure}

\subsection{Extrapolation}
The extrapolation problem of given data is a highly unstable process \cite{Demanet2016}. Recent results, highlight the ability of extended spans when using high arithmetic precision \cite{Bakas2019}. In Figure \ref{fig:extrap}, the highly extended extrapolation span for $f(x)=\sin (x)$ is depicted. The extrapolation errors are starting becoming visible only for $x>73L$. We should highlight, that this is consistent with the corresponding theory as, for this function, the computed ${\mathop{1/r=\lim \sup_{n\to \infty } }}\,\sqrt[n]{|{{a}_{n}}|}$ takes values 0.0178, 0.0169, 0.0161, 0.0152, 0.0145, 0.0137 for the higher values of $n$ (Figure \ref{fig:radius-high}). Accordingly, we may write that $r=1/\text{0}\text{.0137}\simeq \text{72}\text{.99}$, which equals to the observed extrapolation span. Accordingly, the extrapolation lengths for $p=1000$ are 12.141 according to the root test $1/r$ and in the actual computations the errors are $>1$ for $x> 12.150$, and, similarly, for p=500 the root test values is 2.154 and the computed 2.230, as illustrated in Figure \ref{fig:extrap}. Hence, interestingly, utilizing this approach, we may predict not only the behaviour of the approximated unknown function within the given domain, but its extrapolation spans as well, and hence the prediction ability.

\subsection{Numerical Integration}

We calculated the vector $\mathbf{a}$, hence we know an approximation of $f(x)\cong {{a}_{0}}+{{a}_{1}}x+{{a}_{2}}{{x}^{2}}+\cdots +{{a}_{n}}{{x}^{n}}$. By integrating the Taylor polynomial of $f$, the indefinite integral of is \[F(x)\cong {{a}_{0}}x+\frac{{{a}_{1}}{{x}^{2}}}{2}+\frac{{{a}_{2}}{{x}^{3}}}{3}+\cdots +\frac{{{a}_{n}}{{x}^{n+1}}}{n+1}+c.\] The only unknown quantity is $c$, which may be calculated by the supplementary constraint that $F(-L)=0$, hence $c\cong -{{a}_{0}}L-\frac{{{a}_{1}}{{L}^{2}}}{2}-\frac{{{a}_{2}}{{L}^{3}}}{3}-\cdots - \frac{{{a}_{n}}{{L}^{n+1}}}{n+1}$. $f(x)=\sin (x)$, hence $F(x)=-\cos (x)$. Accordingly, ${{F}_{an}}=\int_{-L}^{L}{f(x)dx=-\cos (-L)+\cos (L)}=0$. The proposed scheme offers a direct computation of the integrals, as the vector $\mathbf a$ is known. In Table \ref{tab:table-3}, the vastly low errors of numerical integration are demonstrated, as well as the significance of the studied digits.

\begin{table}[htbp]
\caption{\label{tab:table-3}Numerical integration errors.}
\begin{tabular}{ |p{2cm}||p{1.5cm}|p{1.5cm}|p{1.5cm}|p{1.5cm}|p{1.5cm}|  }
 \hline
 & $p=50$ & $p=100$ & $p=500$ & $p=1000$ & $p=2000$ \\
 \hline
$F_{an}-F_{nu}$ & 1.502e-09 & 3.957e-17 & 1.226e-97 & 2.431e-249 & -1.028e-548\\
\hline
\end{tabular}
\end{table}

\subsection{Numerical Differentiation}
The derivatives of $f$, are inherently computed as\[\mathbf{a}=\left\{ {{a}_{1}},{{a}_{2}},...,{{a}_{n}} \right\}=\left\{ f({{x}_{0}}),\frac{{f}'({{x}_{0}})}{1!},\frac{{f}''({{x}_{0}})}{2!},\cdots ,\frac{{{f}^{(n)}}({{x}_{0}})}{n!} \right\}=\mathbf{df}./\mathbf{n!}\], with $\mathbf{df}$ denoting the vector of the $n$ordinary derivatives of $f$ and $\mathbf{n!}$ the vector of the $n$ factorials. The ${{k}^{th}}<n$ derivative at any other point $x\ne {{x}_{0}}$ may easily be computed By Equation (1), e derive that ${f}'(x)\cong 0+{{a}_{1}}+2{{a}_{2}}x+3{{a}_{3}}{{x}^{2}}+\cdots +n{{a}_{n}}{{x}^{n-1}}$, ${{f}'}'(x)\cong 0+0+2{{a}_{2}}+6{{a}_{3}}x+\cdots +(n-1)n{{a}_{n}}{{x}^{n-2}}$, till
\begin{equation}{{f}^{(k)}}(x)\cong k!{{a}_{k}}{{x}^{k}}+\cdots +\frac{n!}{(n-k)!}{{a}_{n}}{{x}^{n-k}},\end{equation}
where the factors $\left\{ {{a}_{k}},{{a}_{k+1}},...,{{a}_{n}} \right\}$, have already been computed by $\mathbf{a}$. We demonstrate the efficiency of the numerical differentiation in the following example apropos the solution of differential Equations.

\subsection{Solution of Ordinary Differential Equations}
The solution is based on the constitution of the matrices representing the derivatives of $\mathbf{V}$, for example
$\mathbf{dV}=\left[ \begin{matrix}
   0 & 1 & 2{{x}_{1}} & \ldots  & (n-1)x_{1}^{n-2}  \\
   0 & 1 & 2{{x}_{2}} & \ldots  & (n-1)x_{2}^{n-2}  \\
   0 & 1 & 2{{x}_{3}} & \ldots  & (n-1)x_{3}^{n-2}  \\
   \vdots  & \vdots  & \vdots  & \ddots  & \vdots   \\
   0 & 1 & 2{{x}_{n}} & \ldots  & (n-1)x_{n}^{n-2}  \\
\end{matrix} \right]$, and ${{\mathbf{d}}^{2}}\mathbf{V}=\left[ \begin{matrix}
   0 & 0 & 2 & \ldots  & (n-1)(n-2)x_{1}^{n-3}  \\
   0 & 0 & 2 & \ldots  & (n-1)(n-2)x_{2}^{n-3}  \\
   0 & 0 & 2 & \ldots  & (n-1)(n-2)x_{3}^{n-3}  \\
   \vdots  & \vdots  & \vdots  & \ddots  & \vdots   \\
   0 & 0 & 2 & \ldots  & (n-1)(n-2)x_{n}^{n-3}  \\
\end{matrix} \right]$, etc.
By utilizing such matrices, we can easily constitute a system of equations representing the differential equation at points ${{x}_{i}}$. To demonstrate the unified approach for the solution of differential equations, we consider the bending of a simply supported beam \cite{Boresi2011}, with governing equation
\begin{equation}
    EI\frac{{{d}^{4}}w}{d{{x}^{4}}}=q(x)
    \label{eq:Beam}
\end{equation}

where $E$is the modulus of elasticity, $I$the moment of inertia, $w$ the sought solution representing the deflection of the beam, and $q$ the external load. For $E=I=L=1,q(x)=0$, and fixed boundary conditions $w(0)=0,{{\left. \frac{dw}{dx} \right|}_{x=o}}=0,w(L)=1/100,{{\left. \frac{dw}{dx} \right|}_{x=L}}=0$, we may write Equation \ref{eq:Beam} supplemented by the boundary conditions in matrix form by
\[\left[ \begin{matrix}
   0 & 0 & 0 & 0 & 24 & \ldots  & (n-1)(n-2)(n-3)(n-4)x_{1}^{n-5}  \\
   0 & 0 & 0 & 0 & 24 & \ldots  & (n-1)(n-2)(n-3)(n-4)x_{2}^{n-5}  \\
   0 & 0 & 0 & 0 & 24 & \ldots  & (n-1)(n-2)(n-3)(n-4)x_{3}^{n-5}  \\
   \vdots  & \vdots  & \vdots  & \vdots  & \vdots  & \ddots  & \vdots   \\
   0 & 0 & 0 & 0 & 24 & \ldots  & (n-1)(n-2)(n-3)(n-4)x_{n}^{n-5}  \\
   1 & 0 & 0 & 0 & 0 & \cdots  & 0  \\
   0 & 1 & 0 & 0 & 0 & \cdots  & 0  \\
   L & 0 & 0 & 0 & 0 & \cdots  & 0  \\
   0 & L & 0 & 0 & 0 & \cdots  & 0  \\
\end{matrix} \right]\left\{ \begin{matrix}
   {{a}_{0}}  \\
   {{a}_{1}}  \\
   {{a}_{2}}  \\
   \vdots   \\
   {{a}_{n}}  \\
\end{matrix} \right\}=\left\{ \begin{matrix}
   \begin{matrix}
   {{p}_{0}}  \\
   {{p}_{1}}  \\
   {{p}_{2}}  \\
   \vdots   \\
   {{p}_{n}}  \\
\end{matrix}  \\
   {{w}_{0}}  \\
   {{{{w}'}}_{0}}  \\
   {{w}_{L}}  \\
   {{{{w}'}}_{L}}  \\
\end{matrix} \right\}\]

\begin{figure}[ht]   
\centering
\includegraphics[width=0.5\textwidth, keepaspectratio]{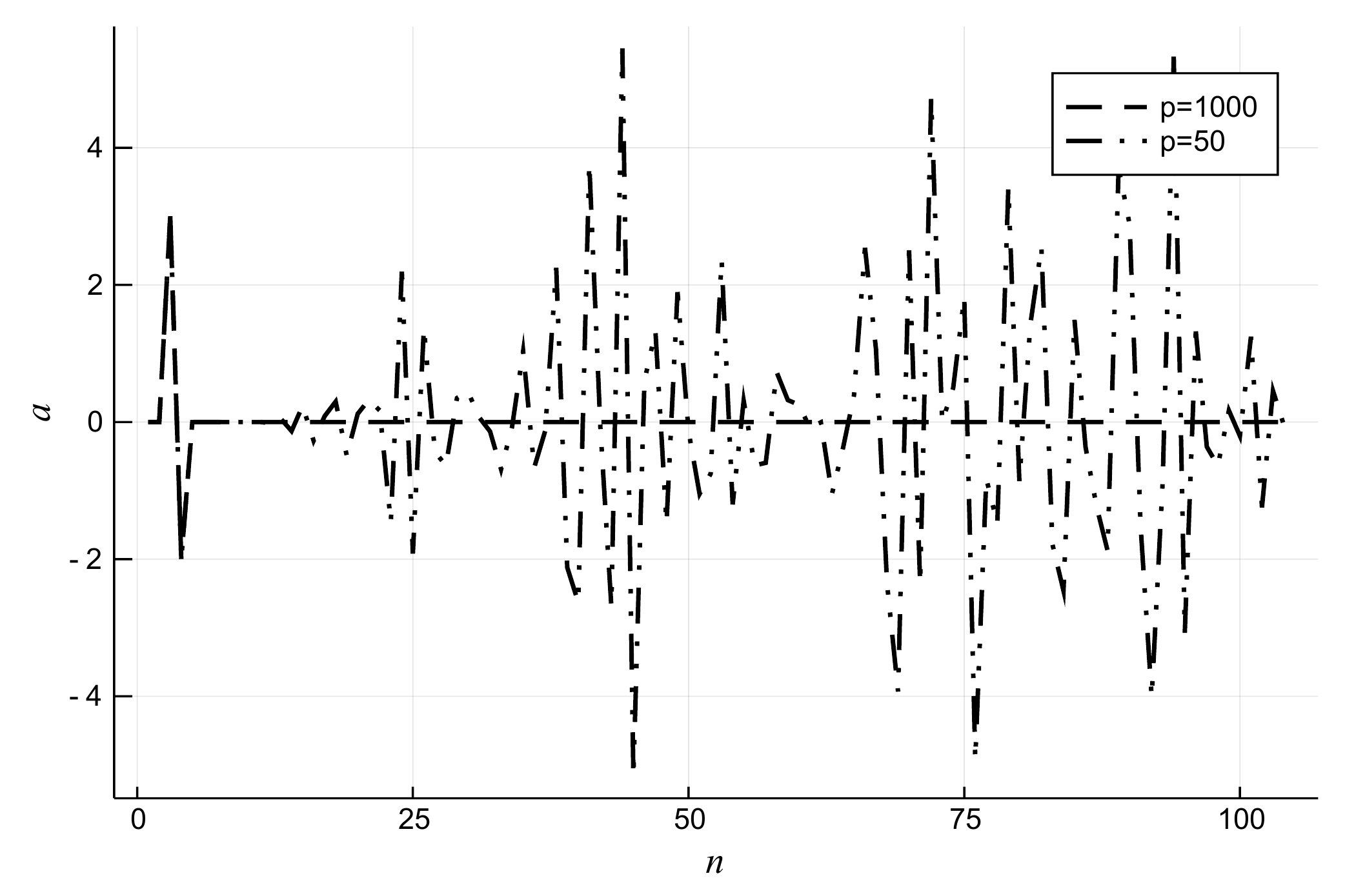}
\caption[]{Calculated $\mathbf{a}$ for $p=50$ and $p=1000$ bits accuracy.}
\label{fig:ODE1}
\end{figure}

Solving for $\mathbf{a}$, and utilizing matrix $\mathbf{V}$, we derive the sought solution by $\mathbf{w}=\mathbf{Va}$.
The exact solution is \[EIw(x)=\frac{-2EI}{{{L}^{3}}}{{x}^{3}}+\frac{3EI}{{{L}^{2}}}{{x}^{2}}\], hence the exact $\mathbf{a}=\left\{ 0,0,3,-2,0,\ldots ,0 \right\}$. In Figure \ref{fig:ODE1}, the ability of high precision ($p=1000$) to identify the exact weights $\mathbf{a}$is revealed, while $p=50$ bits accuracy fails dramatically for such identification. However, they exhibit lower values than the interpolation problem, probably due to the imposition of the boundary conditions.

\subsection{System Identification}

The inverse problems, that is the identification of the system which produced a governing differential law [28], is of great interest as this law describes rigorously the behaviour of a studied system. We demonstrate the ability of high-precision Taylor polynomials for the rapid and precise identification of unknown systems. Let $t$ be an input variable and $s$ a measured response. We may easily compute $\mathbf{a}=\left\{ {{a}_{1}},{{a}_{2}},...,{{a}_{n}} \right\}$, by $\mathbf{a}={{\mathbf{V}}^{-1}}\mathbf{s}$. We assume the existence of a differential operator $T$, such that $T(s)=c$. According to [29], we may write $T$ as a power series by
$T(s)=\sum\limits_{i,j,k=0}^{2}{{{b}_{ijk}}{{s}^{i}}{{{\dot{s}}}^{j}}{{{\ddot{s}}}^{k}}}={{{b}}_{{000}}}{ + }{{{b}}_{{100}}}{s + }{{{b}}_{{010}}}{\dot{s} + }{{{b}}_{{001}}}\ddot{s}{+}{{{b}}_{{200}}}{{{s}}^{2}}{ + }{{{b}}_{{110}}}{s\dot{s}+ }{{{b}}_{{101}}}s\ddot{s}+{{{b}}_{020}}{{{\dot{s}}}^{2}}+{{{b}}_{{011}}}\dot{s}\ddot{s}+{{{b}}_{{002}}}{{\ddot{s}}^{2}}$ and by setting ${{c}'=c-}{{{b}}_{{000}}}$, and assuming a linear approximation, we derive \[1=\frac{{ }{{{b}}_{{100}}}{s + }{{{b}}_{{010}}}{\dot{s} + }{{{b}}_{{001}}}\ddot{s}}{{{c}'}}.\] Applying the later for all ${{x}_{i}}$ and writing the resulting system in matrix form, we obtain 
\begin{equation}[\mathbf{Va}{ + }\mathbf{dVa}{+}{{\mathbf{d}}^{2}}\mathbf{Va}]{{\mathbf{b}}^{T}}=\left\{ \mathbf{1} \right\}\end{equation},
where $\left\{ \mathbf{1} \right\}=\left\{ 1,1,\ldots ,1 \right\}$. Solving for $\mathbf{b}$, we obtain the weights of the derivatives in the differential operator $T(s)$. 

For example if we apply the previous for data of Newton’s second law [30] of motion  $s(t)={{t}^{2}}$, with $s$ indicating space and $t$ time, we may calculate vectors $\mathbf{a}$and solve Equation (4) for $\mathbf{b}$, with ${c}'=1$, and $p=1000$bits precision, we derive that $\mathbf{b}=\left\{ 0,0,1/2 \right\}+O({{10}^{-270}})$, and hence $\frac{1}{2}\ddot{s}=1\to \ddot{s}=2$, which is equivalent with $\ddot{s}=a$, where $a=\frac{F}{m}=2$, which represents the external source which produces $s(t)={{t}^{2}}$.

We assumed that $1={ }{{{b}}_{{100}}}{s + }{{{b}}_{{010}}}{\dot{s} + }{{{b}}_{{001}}}\ddot{s}$, hence by assuming $S=\int{s},SS=\iint{s}$, and integrating in the interval $[0,t]$, we obtain $t+{{c}_{1}}={{{b}}_{{100}}}(S(t)-S(0)) + {{{b}}_{{010}}}{(s(t)-s(0)) + }{{{b}}_{{001}}}(\dot{s}(t)-\dot{s}(0))$, however, $S(0){=s(0)}=\dot{s}(0)=0$. Accordingly, we may write $t={ }{{{b}}_{{100}}}S(t){+ }{{{b}}_{{010}}}{s(t)+ }{{{b}}_{{001}}}\dot{s}(t)$, and if we integrate for a second time in the interval $[0,t]$, we obtain ${}^{{{t}^{2}}}/{}_{2}={ }{{{b}}_{{100}}}(SS(t)-SS(0)){+ }{{{b}}_{{010}}}{(S(t)-S(0))+ }{{{b}}_{{001}}}(s(t)-s(0))$, and because $SS(0)=0$, we have

\begin{equation}
    s(t)={}^{{}^{{{t}^{2}}}/{}_{2}-{{{b}}_{{100}}}SS(t)-{{{b}}_{{010}}}{S(t)}}/{}_{{{{b}}_{{001}}}}{ }
    \label{eq:space1}
\end{equation}

The integrals of $s$, $\int{s}\text{ }$and $\iint{s}$ can be approximated with high accuracy, by utilizing accordingly the procedure discussed in  §3.4, by using the integrals of the obtained Taylor Polynomials \[\int{s}{ }\cong {{a}_{0}}t+\frac{{{a}_{1}}{{t}^{2}}}{2}+\frac{{{a}_{2}}{{t}^{3}}}{3}+\cdots +\frac{{{a}_{n}}{{t}^{n+1}}}{n+1}\]  \[\iint{s}{ }\cong \frac{{{a}_{0}}{{t}^{2}}}{2}+\frac{{{a}_{1}}{{t}^{3}}}{6}+\frac{{{a}_{2}}{{t}^{4}}}{12}+\cdots +\frac{{{a}_{n}}{{t}^{n+2}}}{(n+1)(n+2)}\], as well as the corresponding matrices for all the given ${{t}_{i}}$, \[\mathbf{IV}=\left[ \begin{matrix}
   1 & 1/2 & {{t}_{1}}/3 & \ldots  & t_{1}^{n+1}/(n+1)  \\
   1 & 1/2 & {{t}_{2}}/3 & \ldots  & t_{2}^{n+1}/(n+1)  \\
   1 & 1/2 & {{t}_{3}}/3 & \ldots  & t_{3}^{n+1}/(n+1)  \\
   \vdots  & \vdots  & \vdots  & \ddots  & \vdots   \\
   1 & 1/2 & {{t}_{n}}/3 & \ldots  & t_{n}^{n+1}/(n+1)  \\
\end{matrix} \right]\] \[\mathbf{IIV}=\left[ \begin{matrix}
   1/2 & 1/6 & {{t}_{1}}/12 & \ldots  & t_{1}^{n+2}/(n+1)/(n+2)  \\
   1/2 & 1/6 & {{t}_{2}}/12 & \ldots  & t_{2}^{n+2}/(n+1)/(n+2)  \\
   1/2 & 1/6 & {{t}_{3}}/12 & \ldots  & t_{3}^{n+2}/(n+1)/(n+2)  \\
   \vdots  & \vdots  & \vdots  & \ddots  & \vdots   \\
   1/2 & 1/6 & {{t}_{n}}/12 & \ldots  & t_{n}^{n+2}/(n+1)/(n+2)  \\
\end{matrix} \right]\].

The calculated impact of ${{{b}}_{{001}}}$ for $p=50$ and $p=1000$ bits accuracy is revealed, by the resulting extrapolation curves beyond the observed domain, utilizing Equation \ref{eq:space1}. For $p=50$ bits accuracy, for given data in the domain $[0,1]$ we may extrapolate only up to a short time ($t'=1.343$) after the last given $t_{end}=1.000$, with threshold for errors $<1.000$, while for $p=2000$ bits the corresponding $t'$ attains the remarkably high value of 9.621e+10.

\section{Functions in multiple dimensions}
\label{N-D}

\subsection{Multidimensional Interpolation}

The Taylor series of $f(x,y)$, depending on two variables $x,y\in \Omega $, with $\Omega $ a closed disk about the center ${{x}_{0}},{{y}_{0}}$, may be written utilizing the partial derivatives of $f$ [31], [32], in the form of
$f(x,y)=f(a,b)+(x-a){{f}_{x}}(a,b)+(y-b){{f}_{y}}(a,b)+\frac{1}{2!}({{(x-a)}^{2}}{{f}_{xx}}(a,b)+2(x-a)(y-b){{f}_{xy}}(a,b)+{{(y-b)}^{2}}{{f}_{yy}}(a,b))+\ldots $, which in vector form is written by \[f(\mathbf{x})=f({{\mathbf{x}}_{0}})+{{(\mathbf{x}-{{\mathbf{x}}_{0}})}^{T}}Df({{\mathbf{x}}_{0}})+\frac{1}{2!}{{(\mathbf{x}-{{\mathbf{x}}_{0}})}^{T}}\left\{ {{D}^{2}}f({{\mathbf{x}}_{0}}) \right\}(\mathbf{x}-{{\mathbf{x}}_{0}})+\cdots ,\]
with ${{D}^{2}}~f\,({{\mathbf{x}}_{0}})$, the Hessian matrix at ${{\mathbf{x}}_{0}}$.

Let $n$be the number of given points of $f({{x}_{i}},{{y}_{j}})$, with $i,j\in (1,2,\ldots ,n)$. In order to constitute the approximating polynomial of $f(x,y)$, with high order terms, and formulate the $\mathbf{V}$matrix with dimensions $n\times n$, we consider all possible combinations of $\left\{ {{n}_{i}},{{n}_{j}}\in (0,1,\ldots ,n-1)\mid  {{n}_{i}}+{{n}_{j}}\le n-1 \right\}$. Hence we may write for all the given ${{x}_{i}}$\[\mathbf{V}({{\mathbf{x}}_{i}}\mathbf{,}{{\mathbf{y}}_{j}})=\left[ \begin{matrix}
   1 & {{x}_{1}} & {{y}_{1}} & {{x}_{1}}{{y}_{1}} & x_{1}^{2} & y_{1}^{2} & \ldots  & x_{1}^{{{n}_{k}}}y_{1}^{{{n}_{l}}}  \\
   1 & {{x}_{2}} & {{y}_{2}} & {{x}_{2}}{{y}_{2}} & x_{2}^{2} & y_{2}^{2} & \ldots  & x_{2}^{{{n}_{k}}}y_{2}^{{{n}_{l}}}  \\
   \ldots  & \ldots  & \ldots  & \ldots  & \ldots  & \ldots  & \ddots  & \ldots   \\
   1 & {{x}_{n}} & {{y}_{n}} & {{x}_{n}}{{y}_{n}} & x_{n}^{2} & y_{n}^{2} & \ldots  & x_{n}^{{{n}_{k}}}y_{n}^{{{n}_{l}}}  \\
\end{matrix} \right],\]

with $k+l=n-1$. Thus we can approximate $f$ with $n$ polynomial terms by 	
\begin{equation}
    \mathbf{f=Va}\to \mathbf{a=}{{\mathbf{V}}^{-1}}\mathbf{f}
    \label{eq:mdi1}
\end{equation}
The computation of $\mathbf{a}$ by Equation \ref{eq:mdi1} permits the computation of $f({{\overset{\scriptscriptstyle\frown}{x}}_{i}},{{\overset{\scriptscriptstyle\frown}{y}}_{j}})$, for any ${{\overset{\scriptscriptstyle\frown}{x}}_{i}},{{\overset{\scriptscriptstyle\frown}{y}}_{j}}\in \Omega $, by utilizing the corresponding $\mathbf{\overset{\scriptscriptstyle\frown}{V}}$. 

Let $f(x,y)=\sin (5x)+\cos ({{e}^{2y}})$. We approximate $f$ with $n=300$ random values ${{x}_{i}},{{y}_{i}}\in \left[ -0.5,0.5 \right]$, and later we interpolate $f$ with $n=300$ random values ${{\overset{\scriptscriptstyle\frown}{x}}_{i}},{{\overset{\scriptscriptstyle\frown}{y}}_{j}}_{i}\in \left[ -0.35,0.35 \right]$. In Figure 7, the exact and approximated values $f({{\overset{\scriptscriptstyle\frown}{x}}_{i}},{{\overset{\scriptscriptstyle\frown}{y}}_{j}})$ are depicted, for $p=2000$ and $p=50$ bits accuracy. Apparently, for the same interpolation problem formulation in three dimensions, the computation precision $p$affects dramatically the results. The \[\max \left| f{{({{{\overset{\scriptscriptstyle}{x}}}_{i}},{{{\overset{\scriptscriptstyle}{y}}}_{j}})}_{analytical}}-f{{({{{\overset{\scriptscriptstyle\frown}{x}}}_{i}},{{{\overset{\scriptscriptstyle\frown}{y}}}_{j}})}_{numerical}} \right|\]equals 8.570e-09 for $p=2000$, and 1.286e+01 for $p=50$ bits. The polynomials weight $\mathbf a$ were calculated by firstly computing the $\mathbf{V^{-1}}$ by solving the $\mathbf{V} \backslash \mathbf{I}$, hence $\mathbf{a=V^{-1}f}$, because the $\mathbf{a=V \backslash I}$ exhibited significant errors. The calculation of the inverse of generic  matrices, as well as the solution of systems of Equations in high precision is a topic for future research.

\begin{figure}[ht]   
\centering
\subfloat[$p=2000$ bits]{\includegraphics[width=0.5\textwidth, keepaspectratio]{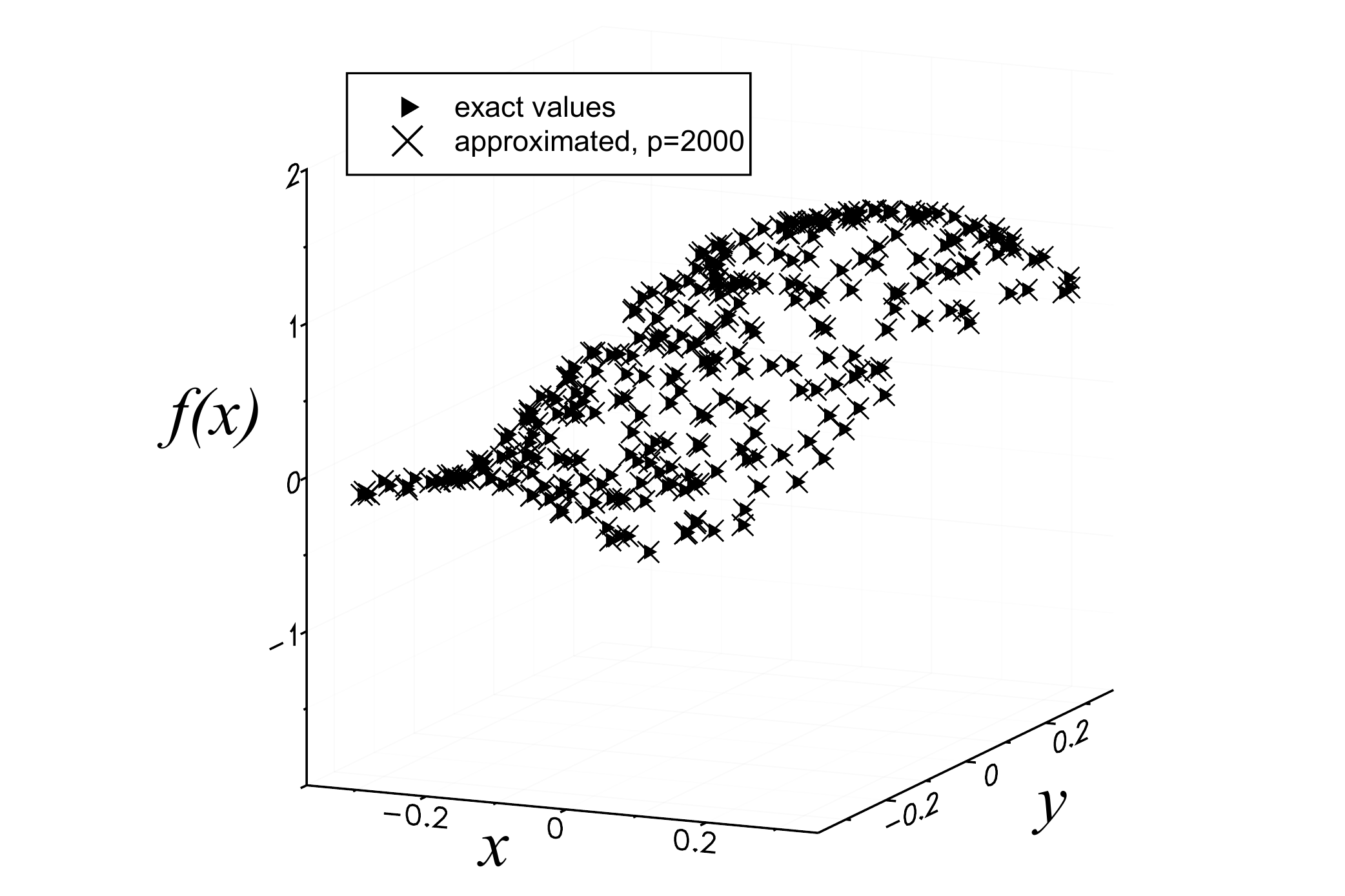}\label{fig:subfig21}}
\subfloat[$p=50$ bits]{\includegraphics[width=0.5\textwidth, keepaspectratio]{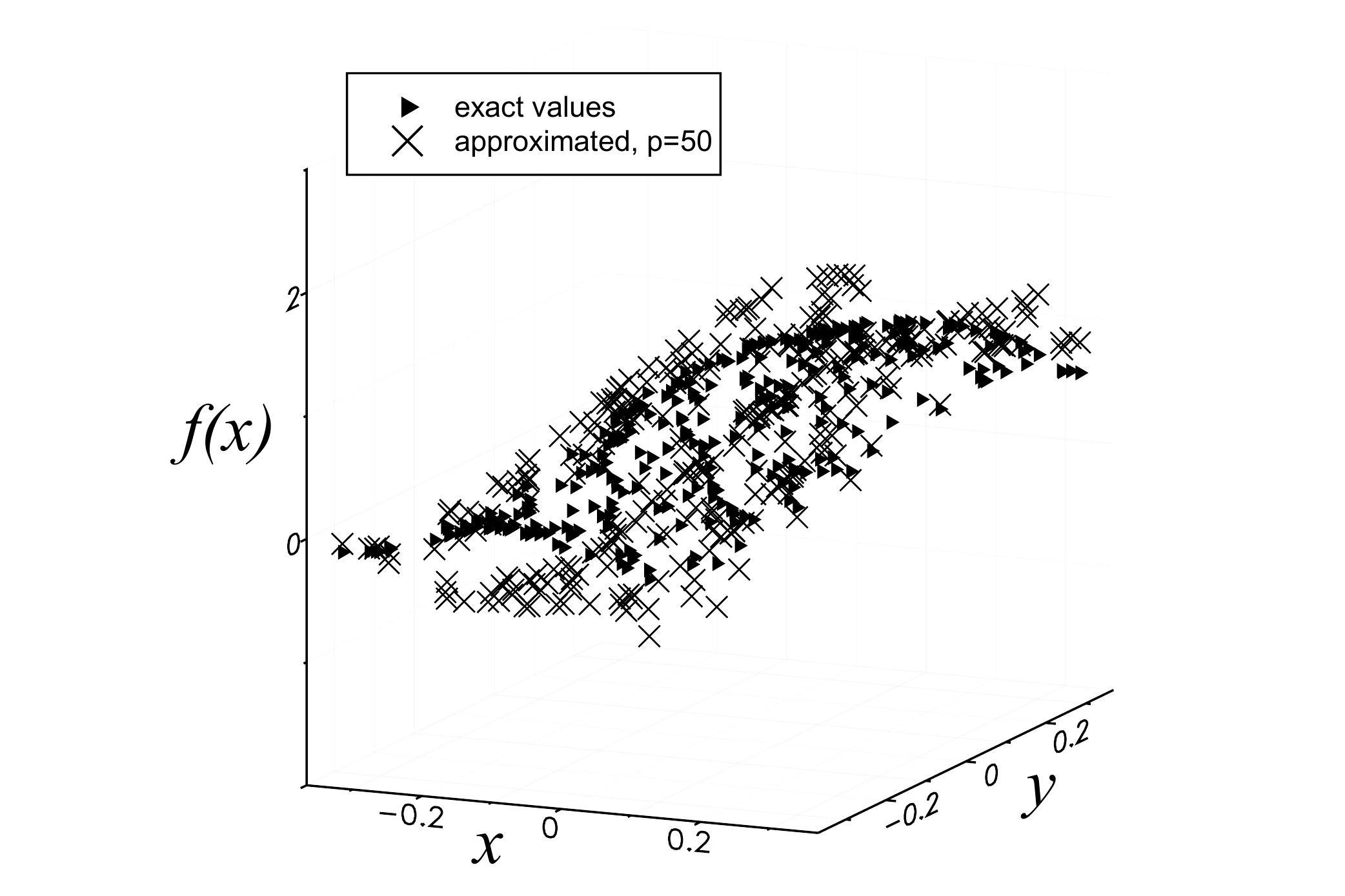}\label{fig:subfig23}}
\caption[Optional caption for list of figures 5-8]{Exact and approximated values of $f$.}
\label{fig:subfigureExample21}
\end{figure}

\subsection{Solution of Partial Differential Equations}

We present the ability of high precision to solve partial differential equations by considering a plate without axial deformations and vertical load $q(x,y)$. The governing equation [33], [34] has the form of
\begin{equation}
    \frac{{{\partial }^{4}}w}{\partial {{x}^{4}}}+2\frac{{{\partial }^{4}}w}{\partial {{x}^{2}}\partial {{y}^{2}}}+\frac{{{\partial }^{4}}w}{\partial {{y}^{4}}}=-\frac{q}{D}
    \label{eq:slab}
\end{equation}
that is ${{\nabla }^{2}}{{\nabla }^{2}}w=-\frac{q}{D}$, with $D:=\frac{2{{h}^{3}}E}{3(1-{{\nu }^{2}})}$, $E$the modulus of elasticity, $v$ the Poisson constant, and $h$ the slab's height.
\par
The sought solution $w(x,y)$ is the slab's deformation within the boundary conditions ${{w}_{b}}({{\mathbf{x}}_{b}},{{\mathbf{y}}_{b}})$ along some boundaries $b=\left\{ 1,2,\ldots  \right\}$. In order to solve Equation \ref{eq:slab}, we approximate \[\mathbf{w=Va}\] using the approximation scheme of Equation \ref{eq:mdi1}, and as the vector $\mathbf{a}$ is constant, we obtain ${{\mathbf{w}}_{{{x}^{4}}}}\mathbf{=}{{\mathbf{V}}_{{{x}^{4}}}}\mathbf{a}$, ${{\mathbf{w}}_{{{y}^{4}}}}\mathbf{=}{{\mathbf{V}}_{{{y}^{4}}}}\mathbf{a}$, ${{\mathbf{w}}_{{{x}^{2}}{{y}^{2}}}}\mathbf{=}{{\mathbf{V}}_{{{x}^{2}}{{y}^{2}}}}\mathbf{a}$, with ${{\mathbf{w}}_{{{x}^{k}}}}_{{{y}^{l}}}$, denoting the partial derivative of $w$, of order $k$over $x$and $l$over $y$, $\frac{{{\partial }^{k+l}}w}{\partial {{x}^{k}}\partial {{y}^{l}}}$, for all given ${{x}_{i}},{{y}_{j}}$ with $i,j\in (1,2,\ldots ,n)$. Utilizing this notation, we may write Equation \ref{eq:slab} for all ${{x}_{i}},{{y}_{j}}$ in matrix form by 
\[\left[ {{\mathbf{V}}_{{{x}^{4}}}}+2{{\mathbf{V}}_{{{x}^{2}}{{y}^{2}}}}+{{\mathbf{V}}_{{{y}^{4}}}} \right]\mathbf{a}=\mathbf{q}.\]

By applying some boundary conditions, we may write for the same $\mathbf{a}$, \par
\[\left[ \begin{matrix}
   {{\mathbf{V}}_{{{x}^{4}}}}+2{{\mathbf{V}}_{{{x}^{2}}{{y}^{2}}}}+{{\mathbf{V}}_{{{y}^{4}}}}  \\
   \mathbf{V}({{x}_{1}},{{y}_{1}})  \\
   {{\mathbf{V}}_{x}}({{x}_{2}},{{y}_{2}})  \\
   \ldots   \\
\end{matrix} \right] \times \mathbf{a}=\left[ \begin{matrix}
   \mathbf{q}  \\
   w({{x}_{1}},{{y}_{1}})  \\
   {{\left. \frac{\partial w}{\partial x} \right|}_{({{x}_{1}},{{y}_{1}})}}  \\
   \ldots   \\
\end{matrix} \right]\to \]

\begin{equation}
    \mathbf{a}={{\left[ \begin{matrix}
   {{\mathbf{V}}_{{{x}^{4}}}}+2{{\mathbf{V}}_{{{x}^{2}}{{y}^{2}}}}+{{\mathbf{V}}_{{{y}^{4}}}}  \\
   \mathbf{V}({{x}_{1}},{{y}_{1}})  \\
   {{\mathbf{V}}_{x}}({{x}_{2}},{{y}_{2}})  \\
   \ldots   \\
\end{matrix} \right]}^{-1}} \times \left[ \begin{matrix}
   \mathbf{q}  \\
   w({{x}_{1}},{{y}_{1}})  \\
   {{\left. \frac{\partial w}{\partial x} \right|}_{({{x}_{1}},{{y}_{1}})}}  \\
   \ldots   \\
\end{matrix} \right].
\label{eq:pdeA}
\end{equation}

By computing $\mathbf{a}$, we then obtain the sought solution as $\mathbf{w=Va}$. \par 

For example, for a simply supported slab, the boundary conditions are $w({{x}_{b}},{{y}_{b}})={{w}_{b}}$ for some boundary $b$. We consider a square slab, with $n=20$divisions per dimension, $dx=1/99$, $L=(n-1)dx$ and $w({{x}_{b}},{{y}_{b}})=0$, at the four linear boundaries. After the computation of $\mathbf{a}$by Equation \ref{eq:pdeA}, we may easily compute the corresponding shear forces, which are defined by 
\[{{Q}_{x}}=-D\frac{\partial }{\partial x}\left( \frac{{{\partial }^{2}}w}{\partial {{x}^{2}}}+\frac{{{\partial }^{2}}w}{\partial {{y}^{2}}} \right), {{Q}_{y}}=-D\frac{\partial }{\partial y}\left( \frac{{{\partial }^{2}}w}{\partial {{x}^{2}}}+\frac{{{\partial }^{2}}w}{\partial {{y}^{2}}} \right).\]

Utilizing the computed $\mathbf{a}$, and matrices ${{\mathbf{V}}_{xxx}},{{\mathbf{V}}_{xyy}},{{\mathbf{V}}_{yxx}},{{\mathbf{V}}_{yyy}}$. Newton equilibrium states that the total shear force at the boundaries should be equal to the total applied force. For constant load over the plate, the Equilibrium errors \[\max \left| \int_{A}{q}(x,y)-\sum{{{Q}_{x,y}}} \right|\]for $p=50$ bits is 6.924e-05 and for $p=2000$ is 2.242 e-591. We observe that there is a big difference, though the errors are small even with $p=50$ bits. Interestingly, utilizing a concentrated, load, by loading the for nodes close to $(0,0)$ the inversion errors $\max \left| {{\left[ \begin{matrix}
   {{\mathbf{V}}_{{{x}^{4}}}}+2{{\mathbf{V}}_{{{x}^{2}}{{y}^{2}}}}+{{\mathbf{V}}_{{{y}^{4}}}}  \\
   \mathbf{V}({{x}_{1}},{{y}_{1}})  \\
   {{\mathbf{V}}_{x}}({{x}_{2}},{{y}_{2}})  \\
   \ldots   \\
\end{matrix} \right]}^{-1}} \times \left[ \begin{matrix}
   {{\mathbf{V}}_{{{x}^{4}}}}+2{{\mathbf{V}}_{{{x}^{2}}{{y}^{2}}}}+{{\mathbf{V}}_{{{y}^{4}}}}  \\
   \mathbf{V}({{x}_{1}},{{y}_{1}})  \\
   {{\mathbf{V}}_{x}}({{x}_{2}},{{y}_{2}})  \\
   \ldots   \\
\end{matrix} \right]-\mathbf{I} \right|$ for $p=50$ bits, is 43.988 and for $p=2000$ is 4.381e-587, further highlighting the significance of accuracy in the calculations.

\section{Conclusions}
\label{sec:conclusions}

System identification and function approximation exist in the core calculations of Physical and Applied Sciences, with implications to other disciplines. Epistemology of scientific discoveries, states that even $1+1=2$ might be falsified \cite{gregory2011arithmetic}. The study of precision in calculations demonstrates illustratively such odd, however fundamental principle. For example, we presented remarkably high extrapolation spans, utilizing a simple representation of the unknown function with Taylor polynomials, by utilizing high arithmetic precision. Approximation errors exhibited great variation in the solutions of Differential Equations, System Identification, and related Numerical Methods. The number of calculation digits are restricted by programming languages' accuracy in bits, however, the utilization of programming structures with extended precision, highlights that certain numerical instabilities stem from the applied computation of the methods’ parameters, and not their theoretical formulation. Interestingly, the approximation errors for the solution of differential equations was even lesser than the interpolation probably due to the imposition of the boundary conditions. We presented the results regarding a variety of numerical methods using function approximation, such as interpolation, extrapolation, numerical differentiation, numerical integration, solution of ordinary and partial differential equations, and system identification, with Taylor polynomials which are in the core foundation of Calculus, as a potential step for the unification of such computational techniques.

\appendix
\section{Programming Code} 
All the results may reproduced by the computer code on GitHub \url{https://github.com/nbakas/TaylorBigF.jl}. The code is in generic form, so as to solve for any numerical problem with the discussed methods.

\bibliographystyle{siamplain}
\bibliography{references}
\end{document}